\documentclass[12pt]{amsart}
\usepackage{amssymb}
\usepackage[all]{xy}

 \rm

\renewcommand{\(}{\left(}
\renewcommand{\)}{\right)}

\renewcommand{\bar}{\overline}

\renewcommand{\phi}{\varphi}

\renewcommand{\span}{{\mathrm{span}}}

\theoremstyle{plain}
\newtheorem{thm}{Theorem}
\newtheorem{lem}[thm]{Lemma}

\theoremstyle{definition}

\theoremstyle{remark}

\title{Raghavan Narsimhan's proof of L.~Schwartz's perturbation theorem}

\author{M.~K.~Vemuri}
\address{Department of Mathematics, West Virginia University,
Morgantown, WV 26506}

\begin{document}

\begin{abstract}
Raghavan Narasimhan outlined a new proof of L.~Schwartz's perturbation
theorem during a course of lectures at IMSc, Chennai in Spring 2007.
The details are given.
\end{abstract}

\maketitle

%\tableofcontents\newpage

Professor Raghavan Narasimhan gave a course of lectures on the
structure of pseudoconvex manifolds at IMSc, Chennai in Spring 2007.
During this course, a new proof of L. Schwartz's perturbation theorem
for operators on Fr\'echet spaces was outlined.  The proof is easy in
the setting of Hilbert spaces, but for Fr\'echet spaces (or even Banach
spaces) the proofs in the literature are quite hard (see
e.g. \cite{NN}).  However, the Fr\'echet space version is the one which
is needed to give easy proofs of the finite dimensionality of various
cohomology groups which occur in complex analysis (in particular, the
easiest proof that there exists a meromorphic function on any compact
Riemann surface).  The proof given here was discovered by Narasimhan
shortly after the second edition of his book \cite{NN} went to press.
Sadly, Professor Narasimhan passed away on October 3, 2015, before a
third edition could be brought out.

There are two main ideas in the proof.  The first is that a compact
set (and hence a compact operator) is small modulo a finite
dimensional subspace.  The second is that a small perturbation of an
onto map is still onto (this is the content of Lemma \ref{L:heart}).
In the setting of Hilbert spaces, this was proved by C.~Neumann in the
19th century using a geometric series argument. Lemma \ref{L:heart} is
a delicate adaptation of his argument, which took over a 100 years to
discover!

\begin{thm}
Let $E, F$ be Fr\'echet spaces, and $f,g:E \to F$ continuous linear maps
such that $g$ is onto, and there exists a neighborhood $U$ of $0$ in
$E$ such that $\bar{f(U)}$ is compact.  Then $f+g$ has closed image of
finite codimension in $F$.
\end{thm}

First of all, observe that it suffices to prove that $(f+g)(E)$ has
finite codimension in $F$, because of the following simple lemma.

\begin{lem}
Let $E,F$ be Fr\'echet spaces and $h:E\to F$ a continuous linear map.
If $h(E)$ has finite codimension in $F$ then $h(E)$ is closed.
\end{lem}

\begin{proof}
Let $F'$ be a complement for $h(E)$ in $F$.  Since $F'$ is finite dimensional,
it is Fr\'echet.  Let $E'=E/\mathrm{Ker}(h)$.  Then $E' \oplus F'$ is Fr\'echet.
Let $\Pi_{F'}: E' \oplus F' \to F'$ denote the projection.  Then $\Pi_{F'}$ is
continuous.  Let $h':E' \to F$ be the map induced by $h$.  Define
$H:E' \oplus F' \to F$ by $H(x,y) = h'(x) + y$.  Then $H$ is continuous,
one-one and onto.  By the open mapping theorem, $H^{-1}$ is continuous.
Moreover, $h(E)=\mathrm{Ker}(\Pi_{F'} \circ H^{-1})$.  Therefore $h(E)$ is
closed. 
\end{proof}

Since $K=\bar{f(U)}$ is compact and $V=g(U)$ is open (by the open mapping
theorem), there exist $y_1, \dots y_n \in K$ such that
$K \subseteq \bigcup_{j=1}^n \( y_j + \frac12 V \)$.  Put
$F'=\span\{y_1, \dots, y_n\}$ and let $f', g':E \to F/F'$ be the induced
maps.  Then $g'$ is onto, $\bar{f'(U)}$ is compact and
$\bar{f'(U)} \subseteq \frac12 g'(U)$.  Thus we are reduced to proving
the following lemma, which is the heart of the matter.

\begin{lem}\label{L:heart}
Let $E, F$ be Fr\'echet spaces, and $f,g: E \to F$ continuous linear maps.
Assume $g$ is onto, and there exists an open symmetric neighborhood $U$
of $0$ in $E$ such that $\bar{f(U)}$ is compact and
$\bar{f(U)} \subseteq \frac12 g(U)$.  Then $h=f+g$ is also onto.
\end{lem}

\begin{proof}
Put $K=\bar{f(U)}$ and $V=g(U)$.  By the open mapping theorem, $V$ is open.
It suffices to show that $h$ is onto $V$.

Let $\{ W_p \}_{p=1}^\infty$ be a fundamental system of neighborhoods of $0$ in
$E$ such that each $W_p$ is open, convex and symmetric.  From the compactness
of $K$, it follows that for each $p$, there exists $n_p$ such that
$K \subseteq \frac12 g(U \cap 2^{n_p} W_p) = g(\frac12 U \cap 2^{n_p - 1} W_p)$.
Discard some of the $W_p$ and reindex them so that
\begin{equation}\label{Ws}
W_{p+1} \subseteq \frac12 W_p
\end{equation}
for all $p$.
The new collection is still a fundamental system of neighborhoods of $0$.

Let $y_0 \in V$.  Then there exists $x_0 \in U$ such that $y_0 = g(x_0)$.
Therefore
\begin{equation*}
y_1 := y_0 - h(x_0) = -f(x_0) \in K = K \cap \frac12 V.
\end{equation*}
Therefore there exists $x_1 \in \frac12 U$ such that
$y_1=g(x_1)$.  Therefore 
\begin{equation*}
y_2 := y_1 - h(x_1) = -f(x_1) \in \frac12 K = \frac12 K \cap \frac14 V.
\end{equation*}
Therefore there exists $x_2 \in \frac14 U$ such that
$y_2=g(x_2)$.  Continuing in this way, we obtain sequences $\{y_j\}_{j=0}^{n_1}$
and $\{x_j\}_{j=0}^{n_1}$ such that
\begin{equation*}
\begin{aligned}
y_{j+1} &= y_j - h(x_j), \\
      y_j &= g(x_j), \\
      y_j &\in 2^{-j+1} K \cap 2^{-j} V, \\
      x_j &\in 2^{-j} U
\end{aligned}
\end{equation*}
Therefore
\begin{equation*}
y_{n_1 + 1} := y_{n_1} - h(x_{n_1}) = -f(x_{n_1}) \in 2^{-n_1} K
=  2^{-n_1} K \cap 2^{-n_1 - 1} V.
\end{equation*}
Therefore there exists $x_{n_1 + 1} \in 2^{-n_1 - 1} U \cap \frac12 W_1$ such that
$g(x_{n_1 + 1}) = y_{n_1 + 1}$.  Continuing in this way, we obtain sequences
$\{y_j\}_{j=n_1 + 1}^{n_2}$ and $\{x_j\}_{j=n_1 + 1}^{n_2}$ such that
\begin{equation*}
\begin{aligned}
y_{j+1} &= y_j - h(x_j), \\
      y_j &= g(x_j), \\
      y_j &\in 2^{-j+1} K \cap 2^{-j} V, \\
      x_j &\in 2^{-j} U \cap 2^{n_1 - j} W_1
\end{aligned}
\end{equation*}

This whole procedure can be further iterated to obtain sequences
$\{y_j\}_{j=0}^\infty$ and $\{x_j\}_{j=0}^{\infty}$ such that
\begin{equation*}
\begin{aligned}
y_{j+1} &= y_j - h(x_j), \\
      y_j &= g(x_j), \\
      y_j &\in 2^{-j+1} K \cap 2^{-j} V, \\
      x_j &\in 2^{-j} U \cap 2^{n_p - j} W_p, \qquad \text{if $j > n_p$},
      \quad p=1, 2, \dots.
\end{aligned}
\end{equation*}

Observe that if $n_p < k < l \le n_{p+1}$ then
\begin{equation*}
\begin{aligned}
x_{k} + \cdots + x_l 
            \in & \;  2^{n_p - k} W_p + \cdots + 2^{n_p - l} W_p \\
\subseteq & \; (2^{n_p - k} + \cdots + 2^{n_p - l}) W_p
                          \qquad\text{(because $W_p$ is convex)} \\
\subseteq & \; W_p.                          
\end{aligned}
\end{equation*}
So if $n_p < k < l \le n_{q}$ then
\begin{equation*}
\begin{aligned}
x_{k} + \cdots + x_l 
            \in & \;  W_p + \cdots + W_{q-1} \\
\subseteq & \;  W_p + \cdots + 2^{p-q} W_p
                        \qquad\text{(by (\ref{Ws}))} \\
\subseteq & \; (1 + \cdots + 2^{p-q}) W_p
                          \qquad\text{(because $W_p$ is convex)} \\
\subseteq & \; 2 W_p.                          
\end{aligned}
\end{equation*}
Therefore $z_j := x_0 + x_1 + \dots + x_j$ is Cauchy.
Let $z=\lim_{j\to\infty} z_j$.  Then
\begin{equation*}
\begin{aligned}
y_0 - h(z)
= & \; \lim_{j\to\infty} (y_0 - h(z_j)) \\
= & \; \lim_{j\to\infty} (y_j - h(x_j)) \\
= & \; \lim_{j\to\infty} y_{j+1}\\
= & \; \lim_{j\to\infty} g(x_{j+1})\\
= & \; 0.
\end{aligned}
\end{equation*}
Therefore $h$ is onto $V$.
\end{proof}

\bibliographystyle{amsplain}
\bibliography{v2-rnpolspt}

\providecommand{\bysame}{\leavevmode\hbox to3em{\hrulefill}\thinspace}
\providecommand{\MR}{\relax\ifhmode\unskip\space\fi MR }
\providecommand{\MRhref}[2]{%
  \href{http://www.ams.org/mathscinet-getitem?mr=#1}{#2}
}
\providecommand{\href}[2]{#2}
\begin{thebibliography}{1}

\bibitem{NN}
Raghavan Narasimhan and Yves Nievergelt, \emph{Complex analysis in one
  variable}, second ed., Birkh\"auser Boston Inc., Boston, MA, 2001.
  \MR{MR1803086 (2002e:30001)}

\end{thebibliography}

\end{document}